\documentclass{tran-l}

 \usepackage[centertags]{amsmath}
 \usepackage{amscd}
 \usepackage{amsfonts}
 \usepackage{amssymb}
 \usepackage{amsthm}
 \usepackage{newlfont}
 \usepackage{graphicx}
 \usepackage[active]{srcltx}

\newcommand{\lc}{\mathcal{L}}
\newcommand{\g}{\mathcal{G}}
\newcommand{\Co}{\operatorname{Cone}}

\theoremstyle{plain}
\newtheorem{thm}{Theorem}[section]

\newtheorem{lem}[thm]{Lemma}
\newtheorem{prop}[thm]{Proposition}
\theoremstyle{definition}
\newtheorem{defn}{Definition}[section]
\theoremstyle{remark}
\newtheorem{rem}{Remark}[section]
\theoremstyle{remark}
\newtheorem{exmp}{Example}[section]
\numberwithin{equation}{section}

\usepackage[cmtip, matrix, arrow]{xy}

\begin{document}

\title{Pre-quantization of Quasi-Hamiltonian Spaces}

\author{ Zohreh Shahbazi }

\address{Department of Mathematics, University of Toronto, Toronto, Ontario,
 Canada}

\email{zohreh@math.utoronto.ca}




\dedicatory{}



\begin{abstract}
This paper develops
  the pre-quantization of Lie group-valued moment maps, and
  establishes its equivalence with the pre-quantization of infinite-dimensional
 Hamiltonian loop group spaces.
\end{abstract}

\maketitle

\section{Introduction}

 \indent The notion of \emph{relative
 gerbes} for smooth maps of manifolds is introduced in ~\cite{ZI}. The equivalence classes of relative gerbes
 are classified by the relative integral cohomology in degree
 three. The differential geometry of relative gerbes, consisting of relative connection, relative connection
curvature, relative Cheeger-Simons differential character, and
relative holonomy are also discussed in ~\cite{ZI}. Inspired by the
pre-quantization of Hamiltonian $G$-manifolds, the main objective of
this paper is to construct a method to pre-quantize the
quasi-Hamiltonian $G$-spaces with group-valued moment maps. For this
purpose, the premise of \emph{relative gerbes} is used. Recall that
a pre-quantization of a symplectic manifold $(M,\omega)$ is a line
bundle L over M with curvature 2-form $\omega$. The symplectic
manifold  $(M,\omega)$ is prequantizable if the 2-form $\omega$ is
integral. In this paper, a notion of a pre-quantization of a space
with G-valued moment map is introduced, and then a similar criterion
for being pre-quantizable is given.
  It is proven that, given two pre-quantizable quasi-Hamiltonian $G$-spaces, their fusion product
  is again pre-quantizable.
Another objective of this paper is to prove that the
pre-quantization of quasi-Hamiltonian G-spaces is equivalent with
the pre-quantization of corresponding infinite-dimensional
 Hamiltonian loop group spaces.

 \indent The organization of this paper is as follows: in Section
2, the relative gerbes are reviewed. In Section 3, an explicit
construction of the basic gerbe for $G=SU(n)$ is given. As well, it
is shown that the construction of the basic gerbe over $SU(n)$ is
equivalent to the construction of the basic gerbe in
Gawedzki-Reis~\cite{MR2003m:81222}. In Section 4, a relative gerbe
for the map $\operatorname{Hol}:\,\mathcal{A}_G(S^1)\rightarrow G$
is constructed, where $\mathcal{A}_G(S^1)$ is the affine space of
connections on the trivial bundle $S^1\times G$.
 In Section
5, quasi-Hamiltonian $G$-spaces ~\cite{MR99k:58062} are reviewed.
The pre-quantization of quasi-Hamiltonian $G$-spaces is introduced
in Section 6, and the pre-quantization conditions for the examples
of
 quasi-Hamiltonian $G$-spaces, described previously, are examined. It is shown that a given
  conjugacy class $\mathcal{C}=G\cdot\exp(\xi)$ of $G$ is pre-quantizable when $\xi\in
  \Lambda^*$, where $\Lambda^*$ is the weight lattice. It is also proven that
  $G^{2h}$ has a unique pre-quantizaion, which enables one to construct
 a finite dimensional pre-quantum line bundle for the moduli space of flat
 connections of a closed oriented surface of genus $h$.

 Further, recall that there is a one-to-one correspondence between quasi-Hamiltonian $G$-spaces and
 infinite dimensional loop group spaces~\cite{MR99k:58062}.
 By extending this correspondence in Section 7, it is proven that pre-quantization of a
 quasi-Hamiltonian $G$-spaces with group-valued moment map coincides with the
 pre-quantization of the corresponding Hamiltonian loop group space.

\section{Review of Relative Gerbes}
\subsection{Gerbes}The main references for this section are
 ~\cite{ZI},~\cite{MR2003f:53086},~\cite{D} and ~\cite{1}.\newline \indent
Let $\mathcal{U}=\{U_i\}_{i\in I}$ be an open cover for a manifold
$M$. It will be convenient to introduce the following notations.
 Suppose that there is a collection of line bundles $L_{i^{(0)},...,i^{(n)}}$ on
 $U_{i^{(0)},...,i^{(n)}}$.
 Consider the inclusion maps,\[\delta_k:\emph{U}_{i^{(0)},...,i^{(n+1)}}\rightarrow
 \emph{U}_{i^{(0)},...,\widehat{i^{(k)}},...,i^{(n+1)}}\quad (k=0,\cdots,n+1)\]
 and define Hermitian line bundles $(\delta L)_{i^{(0)},...,i^{(n+1)}}$
 over $\emph{U}_{i^{(0)},...,i^{(n+1)}}$ by\[\delta L:=\underset{k=0}
 {\overset{n+1}{\bigotimes}}
 (\delta_k^{*}L)^{(-1)^k}.\]Notice that $\delta(\delta L)$ is canonically trivial.
 If one has a unitary section $\lambda_{i^{(0)},...,i^{(n)}}$ of $L_{i^{(0)},...,i^{(n)}}$
 for each $U_{i^{(0)},...,i^{(n)}}\neq\emptyset$, then one can define
$\delta\lambda$ in a similar fashion.
 Note that $\delta(\delta\lambda)=1$ as a section of trivial line bundle.
\begin{defn}
A gerbe on a manifold $M$ on an open cover
$\mathcal{U}=\{\emph{U}_{i}\}_{i\in I}$ of $M$ is defined by
Hermitian line bundles $L_{ii'}$ on each $\emph{U}_{ii'}$ such that
$L_{ii'}\cong L_{i'i}^{-1}$, and a unitary section $
\theta_{ii'i''}$ of $\delta L$ on $\emph{U}_{ii'i''}$ such that
$\delta\theta=1$ on $\emph{U}_{ii'i''i'''}$. Denote this data as
$\mathcal{G}=(\mathcal{U},L,\theta)$.\end{defn} Denote the set of
all gerbes on $M$ as $Ger(M)$. $Ger(M)$ has a natural group
structure and is classified with
$H^3(M,\mathbb{Z})$.\goodbreak\begin{defn} A quasi-line bundle for
the gerbe $\mathcal{G}$ on a manifold $M$ on the open cover
$\mathcal{U}=\{U_{i}\}_{i\in I}$ is defined
as:\begin{enumerate}\item a Hermitian line bundle $E_{i}$ over each
$U_{i}$;
\item Unitary sections $\psi_{ii'}$ of
\[(\delta E^{-1})_{ii'}\otimes L_{ii'}\]such that
$\delta\psi=\theta$.\end{enumerate} Denote this quasi-line bundle as
$\mathcal{L}=(E,\psi).$
\end{defn} Any two quasi-line bundles over a given gerbe differ by
a line bundle.
\begin{defn} A relative gerbe for $\Phi\in C^{\infty}(M,N)$ is a
pair $(\mathcal{L},\mathcal{G})$ where $\g$ is a gerbe over $N$ and
$\mathcal{L}$ is a quasi-line bundle for $\Phi^*\g$. Denote the set
of all relative gerbes for the map $\Phi$ as $Ger(\Phi)$.\end{defn}
$Ger(\Phi)$ has a group structure and is classified with the third
degree relative integral cohomology of the map $\Phi$, i.e.,
$$Ger(\Phi)\cong H^3(\Phi,\mathbb{Z}).$$
\section{Gerbes over a Compact Lie Group }It is
well-known that for a compact, simple, simply connected Lie group
the integral cohomology $H^{\bullet}(G,\mathbb{Z})$ is trivial in
degree less than three, while $H^3(G,\mathbb{Z})$ is canonically
isomorphic to $\mathbb{Z}$. The gerbe corresponding to the generator
of $H^3(G,\mathbb{Z})$ is called the basic gerbe over $G$. In this
section, an explicit construction of the basic gerbe for $G=SU(n)$
is given. This gerbe plays an important role in pre-quantization of
the quasi-Hamiltonian $G$-spaces.\subsection{Some Notations from Lie
Groups}Let $G$ be a compact, simple simply connected Lie group, with
a maximal torus $T$. Let $\mathfrak{g}$ and $\mathfrak{t}$ denote
the Lie algebras of $G$ and $T$, respectively. Denote by
$\Lambda\subset \mathfrak{t}$ the integral lattice, given as the
kernel of
$$\exp:\mathfrak{t}\rightarrow T.$$
Let $\Lambda^*=Hom(\Lambda,\mathbb{Z})\subset \mathfrak{t}^*$ be its
dual weight lattice. Recall that any $\mu\in \Lambda^*$ defines a
homomorphism
$$h_{\mu}:T\rightarrow U(1),\,\exp\xi\mapsto e^{2\pi\sqrt{-1}
\langle\mu,\xi\rangle}.$$ This identifies $\Lambda^*=Hom(T,U(1))$.
Let $\mathcal{R}\subset \Lambda^*$ be the set of roots, i.e., the
non-zero weights for the adjoint representation. Define
$$\mathfrak{t}^{reg}:=\mathfrak{t} \setminus\bigcup_{\alpha\in \mathcal{R}}
ker\alpha.$$ The closures of the connected components of
$\mathfrak{t}^{reg}$ are called Weyl chambers. Fix a Weyl chamber
$\mathfrak{t}_+$. Let $\mathcal{R}_+\subset \Lambda^*$ be the set of
the positive roots, i.e., roots that are non-negative on
$\mathfrak{t}_+$. Then $\mathcal{R}=\mathcal{R}_+\cup
-\mathcal{R}_+$. A positive root is called simple, if it cannot be
written as the sum of positive roots. Denote the set of simple roots
by $\mathcal{S}$.  The set of simple roots $\mathcal{S}\subset
\mathcal{R}_+$ forms a basis of $\mathfrak{t}$, and
$$\mathfrak{t}_{+}=\{\xi\in\mathfrak{t}\mid\langle\alpha,\xi\rangle\geq 0,
\forall\alpha\in S\}.$$ Any root $\alpha\in \mathcal{R}$ can be
uniquely written as
$$\alpha=\Sigma k_i\alpha_i,\quad k_i\in \mathbb{Z},\alpha_i\in\mathcal{S}.$$
The hight of $\alpha$ is defined by
$\operatorname{ht}(\alpha)=\Sigma k_i$. Since $\mathfrak{g}$ is
simple, there is a unique root $\alpha_0$ with
$\operatorname{ht}(\alpha)\geq \operatorname{ht}(\alpha_0)$ for all
$\alpha\in \mathcal{R}$, which is called the lowest root. The
fundamental alcove is defined as
$$\mathfrak{A}=\{\xi\in\mathfrak{t}_{+}\mid\langle\alpha_0,\xi\rangle\geq-1\}.$$
 The basic inner product
on $\mathfrak{g}$ is the unique invariant inner product such that
$\alpha.\alpha=2$ for all long roots $\alpha$, which is used here to
identify $\mathfrak{g}^*\cong\mathfrak{g}$. The mapping
$\xi\rightarrow \operatorname{Ad} G(\exp\xi)$ is a homeomorphism
from $\mathfrak{A}$ onto $G/\operatorname{Ad} G$, the space of the
cojugacy classes in $G$. Therefore, the fundamental alcove
parameterizes conjugacy classes in $G$~\cite{DK}. Denote the
quotient map by $q:G\rightarrow \mathfrak{A}$.\newline \indent Let
$\theta^L,\theta^R\in \Omega^1(G,\mathfrak{g})$ be the left and
right invariant Maurer Cartan forms. If $L_g$ and $R_g$ denote left
and right multiplication by $g\in G$, then the values of
$\theta^L_{g}$ and $\theta^R_{g}$ at $g$ are given by
\[\theta^L_{g}=dL_{g^{-1}} : TG_g\rightarrow
TG_e\cong\mathfrak{g},\quad \theta^R_{g}=dR_{g^{-1}} :
TG_g\rightarrow TG_e\cong\mathfrak{g}.\]For any $g\in G$,
\[\theta^L_{g}=Ad_g(\theta^R_{g}).\]For any invariant inner product
$B$ on $\mathfrak{g}$, the form
\begin{equation}\label{3-form}\eta:=\frac{1}{12} B(\theta^L,
[\theta^L,\theta^L])\in \Omega^3(G)\end{equation} is bi-invariant
since the inner product is invariant. Any bi-invariant form on a Lie
group is closed. Therefore, $\eta$ is a closed 3-form and its
cohomology class represents the generator of
$H^3(G,\mathbb{R})=\mathbb{R}$ if one assumes that $G$ is compact
and simple. If, in addition, $G$ is simply connected, then
$H^3(G,\mathbb{Z})=\mathbb{Z}$, and one can normalize the inner
product such that $[\eta]$ represents an integral generator
~\cite{MR2001i:53140}. $B$ is called the inner product at level
$\lambda>0$ if $B(\xi,\xi)=2\lambda$ for all short lattice vectors
$\xi\in \Lambda$. The inner product at level $\lambda=1$ is called
the basic inner product. (It is related to the Killing form by a
factor $2c_{\mathfrak{g}}$, where $c_{\mathfrak{g}}$ is the dual
Coxeter number of $\mathfrak{g}$.) Suppose that $G$ is simply
connected and simple. It is known that the 3-form defined by
Equation \ref{3-form} is integral if and only if its level $\lambda$
of $B$ is an integer.
\subsection{Standard Open Cover of G} Let $\mu_0,\cdots,\mu_d$ be
the vertices of $\mathfrak{A}$, with $\mu_0=0$. Let $\mathfrak{A}_j$
be the complement of the closed face opposite to the vertex $\mu_j$.
The standard open cover of $G$ is defined by the pre-images
$V_j=q^{-1}(\mathfrak{A}_j)$. Denote the centralizer of $\exp\mu_j$
by $G_j$. Then, the flow-out $S_j=G_j.\exp(\mathfrak{A}_j)$ is an
open subset of $G_j$, and it is a slice for the conjugation action
of $G$. Therefore,
$$G\times_{G_{j}}S_j=V_j.$$More generally, let
$\mathfrak{A}_I=\cap_{j\in I}\mathfrak{A}_j$ and
$V_I=q^{-1}(\mathfrak{A}_I)$. Then, $S_I=G_I.\exp(\mathfrak{A}_I)$
is a slice for the conjugation action of $G$, and therefore
$$G\times_{G_{I}}S_I=V_I.$$ Denote the projection to the base
by $$\pi_I: V_I\rightarrow G/G_I.$$\begin{lem}$\eta_G$ is exact over
each of the open subsets $V_j$.\end{lem}\begin{proof} $S_j':=
G_j\cdot(\mathfrak{A}_j-\mu_j)$ is a star-shaped open neighborhood
of 0 in $\mathfrak{g}_j$ and is $G_j$-equivariantly diffeomorphic
with $S_j$. One can extend this retraction from $S_j$ onto
$\exp(\mu_j)$ to a $G$-equivariant retraction from $V_j$ onto
$\mathcal{C}_j=q^{-1}(\mu_j)$. But, since
$d_G\omega_{\mathcal{C}_j}+\iota^*_{\mathcal{C}_j}\eta_G=0$, then
$\eta_G$ is exact over $V_j$.
\end{proof}
Let $\iota_j :\mathcal{C}_j\rightarrow V_j$ and $\pi_j
:V_j\rightarrow G/G_j=\mathcal{C}_j$ denote the inclusion and the
projection respectively. The retraction from $V_j$ onto
$\mathcal{C}_j$ defines a $G$-equivariant homotopy operator
$$h_j :\Omega^p(V_j)\rightarrow \Omega^{p-1}(V_j).$$ Thus,
$$d_Gh_j+h_jd_G=Id-\pi_j^*\iota_j^*.$$Define the equivariant 2-form $\varpi_j$ on $V_j$ by
$(\varpi_j)_G=h_j\eta_G-\pi_j^*\omega_{\mathcal{C}_j}$. Write
$(\varpi_j)_G=\varpi_j-\theta_j$, where $\varpi_j\in \Omega^2(V_j)$
and $\theta_j\in \Omega^0(V_j,\mathfrak{g})$. For any conjugacy
class $\mathcal{C}\subset V_j$,
$\iota^*(\varpi_j)_G+\omega_{\mathcal{C}}$ is an equivariantly
closed 2-form with $\theta_j$ as its moment map. Therefore,
$\iota^*(\varpi_j)_G+\omega_{\mathcal{C}}=\theta_j^*(\omega_{\mathcal{O}})_G$,
where $(\omega_{\mathcal{O}})_G$ is the symplectic form on the
(co)-adjoint orbit
$\mathcal{O}=\theta_j(\mathcal{C})$.\begin{prop}Over $V_{ij}=V_i\cap
V_j$, $\theta_i-\theta_j$ takes values in the adjoint orbit
$\mathcal{O}_{ij}$ through $\mu_i-\mu_j$. Furthermore,
$$(\varpi_i)_G-(\varpi_j)_G=\theta_{ij}^*(\omega_{\mathcal{O}_{ij}})_G$$where
$\theta_{ij} :=\theta_i-\theta_j : V_{ij}\rightarrow
\mathcal{O}_{ij}$, and $(\omega_{\mathcal{O}_{ij}})_G$ is the
equivariant symplectic form on the orbit.\end{prop}\begin{proof} Let
$\nu : \mathfrak{A}_j\rightarrow \mathfrak{t}$ be the inclusion map.
Then, $$\widetilde{h}_j\circ
(\exp\mid_{\mathfrak{A}_j})\frac{1}{2}(\theta^L+\theta^R)=\widetilde{h}_j\circ
d\nu=\nu-\mu_j$$ where $\widetilde{h}_j$ is the homotopy operator
for the linear retraction of $\mathfrak{t}$ onto $\mu_j$. This
proves that $(\exp\mid_{\mathfrak{A}_j})^*\theta_j=\nu-\mu_j$.
Therefore, for $\xi\in \mathfrak{A}_{ij}$,
$$\theta_{ij}(\exp\xi)=(\xi-\mu_i)-(\xi-\mu_j)=\mu_j-\mu_i.$$Therefore,
$\theta_{ij}$ takes values in the adjoint orbit through
$\mu_j-\mu_i$ by equivariance. The difference $\varpi_i-\varpi_j$
vanishes on $T$, and is, therefore, determined by its contractions
with generating vector fields. But, $\theta_{ij}$ is a moment map
for $\varpi_i-\varpi_j$, hence $\varpi_i-\varpi_j$ equals to the
pull-back of the symplectic form on
$\mathcal{O}_{ij}$.\end{proof}\subsection{Construction of the Basic
Gerbe}\label{basic gerbe}Let $G$ be a compact, simple, simply
connected Lie group, and $B$ be an invariant inner product
 at integral level $k>0$. Use $B$ to identify
 $\mathfrak{g}\cong \mathfrak{g}^*$ and  $\mathfrak{t}\cong
 \mathfrak{t}^*$.
 Assume that under this identification, all vertices of the alcove are contained in
the weight lattice $\Lambda^*\subset\mathfrak{t}$. This is automatic
if $G$ is the special unitary group $A_d=SU(d+1)$ or the compact
symplectic group $C_d=Sp(2d)$. In general, the following table lists
the smallest integer $k$ with this
property~\cite{MR39:1590}: \\[5.pt]
\begin{center}\begin{tabular}{c|c|c|c|c|c|c|c|c|c}
$G$& $A_d$& $B_d$& $C_d$& $D_d$& $E_6$& $E_7$& $E_8$& $F_4$& $G_2$\\
\hline
$k$& 1& 2& 1& 2& 3& 12& 60& 6& 2\\
\end{tabular}
\\[40.pt]\end{center}
For constructing the basic gerbe over $G$, pick the standard open
cover of $G$, $\mathcal{V}=\{V_i, i=0,\cdots,d\}$. For any $\mu\in
\Lambda^*$, with stabilizer $G_{\mu}$, define a line bundle
$$L_{\mu}=G\times_{G_{\mu}} \mathbb{C}_{\mu}$$ with the unique
left invariant connection $\nabla$, where $\mathbb{C}_{\mu}$ is the
1-dimensional $G_{\mu}$-representation with infinitesimal character
$\mu$. $L_{\mu}$ is a $G$-equivariant pre-quantum line bundle for
the orbit $\mathcal{O}=G\cdot\mu$.
Therefore,$$\frac{i}{2\pi}curve_G(\nabla)=(\omega_{\mathcal{C}})_G
:=\omega_{\mathcal{O}}-\Phi_{\mathcal{O}}$$where
$\omega_{\mathcal{O}}$ is a symplectic form for the inclusion map
$\Phi_{\mathcal{O}} :\mathcal{O}\rightarrow \mathfrak{g}^*$. Define
line bundles $$L_{ij} :=\theta_{ij}^*(L_{\mu_j-\mu_i})$$equipped
with the pull-back connection. In three fold intersection $V_{ijk}$,
the tensor product $(\delta L)_{ijk}=L_{jk}L_{ik}^{-1}L_{ij}$ is the
pull-back of the line bundle over $G/G_{ijk}$, which is defined by
the zero weight $$(\mu_k-\mu_j)-(\mu_k-\mu_i)+(\mu_j-\mu_i)=0$$of
$G_{ijk}$. Therefore, it is canonically trivial with trivial
connection. The trivial sections $t_{ijk}=1$ satisfy $\delta t=1$
and $(\delta\nabla)t=0$. Define $(F_j)_G=(\varpi_j)_G$.
Since$$\delta(F)_G=\theta_{ij}^*(\omega_{\mathcal{O}_{ij}})_G=\frac{1}{(2\pi
\sqrt{-1})}curve_G(\nabla^{ij}),$$ then $\g=(\mathcal{V},L,t)$ is a
gerbe with connection $(\nabla,\varpi)$. The construction of the
basic gerbe is discussed in more general cases
in~\cite{EM,MR1968268}.
 \subsection{The
Basic Gerbe Over SU(n)} This Section, shows that the construction of
the basic gerbe over $SU(n)$, as discussed in the previous section,
is equivalent to the construction of the basic gerbe in
Gawedzki-Reis~\cite{MR2003m:81222}.\newline \indent The
\emph{special unitary group} is the classical group:$$SU(n)=\{A\in
U(n)\mid \det A=1\},$$which is a compact connected Lie group of
dimension equal to $n^2-1$ with Lie algebra equal to the
space:$$\mathfrak{su}(n)=\{A\in
L_{\mathbb{C}}(\mathbb{C}^n,\mathbb{C}^n)\mid A^{\ast}+A=0
\,\mbox{and}\, tr A=0\}.$$Any matrix $A\in SU(n)$ is conjugate to a
diagonal matrix with entries
$$\operatorname{diag}\big(\exp((2\pi
\sqrt{-1})\lambda_1(A)),\cdots, \exp((2\pi
\sqrt{-1})\lambda_n(A))\big)$$ where
$\lambda_1(A),\cdots,\lambda_n(A)\in \mathbb{R}$ are normalized by
the identity $\Sigma_{i=1}^{n}\lambda_i(A)=0$, and
\begin{equation}\lambda_1(A)\geq \lambda_2(A)\geq \cdots\geq
\lambda_n(A)\geq\lambda_1(A)-1.\end{equation}Consider the following
maximal torus of $SU(n)$$$T=\{A\in SU(n)\mid \mbox{$A$ is
diagonal}\}.$$Let $\mathfrak{t}$ be the Lie algebra of $T$. Thus,
 $\mathfrak{t}\cong\{\lambda\in\mathbb{R}^n\mid\Sigma_{i=1}^n\lambda_i=0\}.$
The roots $\alpha\in \mathcal{R}\subset \mathfrak{t}^*$ are the
linear maps:$$\alpha_{ij}:\mathfrak{t}\rightarrow
\mathbb{R},\,(\lambda_1,\cdots,\lambda_n)\mapsto\lambda_i-\lambda_j,\,i\neq
j,$$and the set of simple roots is
$$\mathcal{S}=\{\alpha_{1,2},\alpha_{2,3},\cdots,\alpha_{n-1,n}\}.$$The
lowest root is $\alpha_{n,1}$(~\cite{MR2003c:22001}, Appendix C).
Choose the following Weyl
chamber$$\mathfrak{t}_+=\{(\lambda_1,\cdots,\lambda_n)\in
\mathfrak{t}\mid \lambda_1\geq\lambda_2\geq\cdots\geq\lambda_n\}.$$
In the above case the fundamental alcove is
$$\mathfrak{A}=\{(\lambda_1,\cdots,\lambda_n)\in \mathfrak{t}\mid
\lambda_1\geq\lambda_2\geq\cdots\geq\lambda_n\geq\lambda_1-1\}.$$The
basic inner product on $\mathfrak{t}$ is induced from the standard
basic inner product on $\mathbb{R}^n$. One can use this inner
product to identify $\mathfrak{t}\cong\mathfrak{t}^*$. Under this
identification $\alpha_{i,j}=e_i-e_j$, where $\{e_i\}_{i=1}^n$ is
the standard basis for $\mathbb{R}^n$. The fundamental weights are
given by $$\mu_i=\{\lambda\in \mathfrak{A}\mid
\lambda_1=\lambda_2=\cdots=\lambda_i>\lambda_{i+1}=\cdots
=\lambda_n=\lambda_1-1\}.$$
$$SU(n)^{reg}=\{A\in SU(n)\mid \mbox{all eigenvalues of $A$ have multiplicity one}\}\cong G
\times_T \,\operatorname{int}\mathfrak{A}\cong G/T.$$ For $i\in
\{1,\cdots,n\}$, define
$$\mathfrak{A}_i:= \{\lambda\in \mathfrak{A}\mid
\lambda_1\geq\cdots \geq\lambda_i>\lambda_{i+1}\geq\cdots
\geq\lambda_n\geq\lambda_1-1\}.$$Thus, the standard open cover for
$SU(n)$ is $\mathcal{V}=\{V_i\}_{i=1}^n$, where
$V_i=q^{-1}(\mathfrak{A}_i)$. Each $SU(n)_{ij}$ is isomorphic to
$U(n-1)$ with the center isomorphic to $U(1)$. Over the set of
regular elements all the inequalities are strict, and one has n
equivariant line bundles $E_1,\cdots,E_n$ defined by the eigenlines
for the eigenvalues $\exp((2\pi \sqrt{-1})\lambda_i(A))$. For $i<j$,
the tensor product $E_{i+1}\otimes\cdots \otimes E_j\rightarrow
SU(n)^{reg}$ extends to a line bundle $E_{ij}\rightarrow V_{ij}$.
For $i<j<k$, one can have a canonical isomorphism $E_{ij}\otimes
E_{jk}\cong E_{ik}$ over $V_{ijk}$. These line bundles together with
corresponding isomorphisms define a gerbe over $SU(n)$, in
Gawedzki-Reis sense, which represents the generator of
$H^3(SU(n),\mathbb{Z})$. Each $E_i$ is equal to $G\times_T
\mathbb{C}_{\nu_i}$ for some $\nu\in \Lambda^*$. In fact, by using
the standard action of $T\subset SU(n)$ on $\mathbb{C}^n$, one can
see that
$$\nu_i=e_i-\frac{1}{n}(1,\cdots,1).$$ Since $\mu_i=\Sigma_{k=1}^i e_k-\frac{i}{n}
\Sigma_{k=1}^n e_k,$ $\mu_i=\Sigma_{k=1}^i \nu_k.$ Recall from
Section \ref{basic gerbe} that $L_{ij}
:=\theta_{ij}^*(L_{\mu_j-\mu_i})$ on $V_{ij}$. Thus, for $i<j$
\begin{eqnarray*}L_{ij}
&=&\theta_{ij}^*(L_{\mu_j-\mu_i})\\&=&\theta_{ij}^*(L_{\Sigma_{k=1}^j
\nu_k-\Sigma_{k=1}^i \nu_k})\\&=&\theta_{ij}^*(L_{\Sigma_{k=i+1}^j
\nu_k})=E_{ij}\end{eqnarray*} \section{The Relative Gerbe for
$\operatorname{Hol}:\label{rhol} \mathcal{A}_G(S^1)\rightarrow
G$}Denote the affine space of connection on the trivial bundle
$S^1\times G$ by $\mathcal{A}_G(S^1)$. Thus,
$\mathcal{A}_G(S^1)=\Omega^1(S^1,\mathfrak{g})$. The loop group
$LG=Map(S^1,G)$ acts on $\mathcal{A}_G(S^1)$ by gauge
transformations:\begin{equation}\label{hol}g\cdot
A=Ad_g(A)-g^*\theta^R.\end{equation}Taking the holonomy of a
connection defines a smooth map $$\operatorname{Hol}:
\mathcal{A}_G(S^1)\rightarrow G$$with equivariance property
$\operatorname{Hol}(g\cdot A)=Ad_{g(0)}\operatorname{Hol}(A)$. If
$\mathfrak{g}$ carries an invariant inner product $B$, write
$L\mathfrak{g}^*$ instead of $\mathcal{A}_G(S^1)$ using the natural
pairing between $\Omega^1(S^1,\mathfrak{g})$ and
$L\mathfrak{g}=\Omega^0(S^1,\mathfrak{g})$. Let's refer to this
action as the coadjoint action. However, notice that the action
\ref{hol} is not the point-wise action. Recall from Section
\ref{basic gerbe},
\\ \indent a) An open cover $\mathcal{V}=\{V_0,\cdots,V_d\}$ of
$G$ such that $V_j/G=\mathfrak{A}_j$, where $d=rank G$.\\
\indent b) For each $V_j$, a unique $G$-equivariant deformation
retraction on to a conjugacy class
$\mathcal{C}_j=G\cdot\exp(\mu_j)$, where $\mu_j$ is the vertex of
$\mathfrak{A}_j$. This deformation retraction descends to the linear
retraction of $\mathfrak{A}_j$ to $\mu_j$.\\ \indent c) 2-forms
$\varpi_j\in \Omega^2(V_j)$, with $d\varpi_j=\eta\mid_{V_j}$, such
that the pull-back onto $\mathcal{C}_j$ is the invariant 2-form for
the conjugacy class $\mathcal{C}_j$.\begin{lem}There exists a unique
$LG$-equivariant retraction from
$\widetilde{V}_j:=\operatorname{Hol}^{-1}(V_j)$ onto the coadjoint
orbit $\mathcal{O}_j=LG\cdot\mu_j$, descending to the linear
retraction of $\mathfrak{A}_j$ onto
$\mu_j$.\end{lem}\begin{proof}The holonomy map sets up a one-to-one
correspondence between the sets of $G$-conjugacy classes and
coadjoint $LG$-orbits, hence both are parameterized by points in the
alcove. The evaluation map $LG\rightarrow G,\,g\mapsto g(1)$
restricts to an isomorphism $(LG)_j\cong G_j$ ~\cite{MR2002g:53151}.
Hence,$$LG\times_{(LG)_j}\widetilde{S_j}=\widetilde{V_j}$$where
$\widetilde{S_j}=\operatorname{Hol}^{-1}(S_j)$ and
$S_j=G_j\cdot\exp(\mathfrak{A}_j)$. Therefore, the unique
equivariant retraction from $V_j$ onto $\mathcal{C}_j$, which
descends to the linear retraction of $\mathfrak{A}_j$ onto the
vertex $\mu_j$, lifts to the desired retraction.
\end{proof}Consider the following commutative diagram:\\*\[\begin{CD}
   \widetilde{V}_j  @>\widetilde{\pi}_j>>     \mathcal{O}_j\\
   @V\operatorname{Hol}VV   @VV\operatorname{Hol}V\\
   V_j@>\pi_j>>\mathcal{C}_j
   \end{CD}\]\\*where $\widetilde{\pi}_j:\widetilde{V}_j
   \rightarrow \mathcal{O}_j$ is the projection that is obtained from
   retraction. Let $\sigma_j\in \Omega^2(\widetilde{V}_j)$ denote
   the pull-backs under $\widetilde{\pi}_j$ of the symplectic forms
   on $\mathcal{O}_j$.\begin{lem}On overlaps
   $\widetilde{V}_j\cap\widetilde{V}_{j'}$,
   $\sigma_j-\sigma_{j'}=\operatorname{Hol}^*(\varpi_j-\varpi_{j'})$.\end{lem}
   \begin{proof}Both sides are closed $LG$-invariant forms, for
   which the pull-back to $\mathfrak{t}\subset L\mathfrak{g}^*$
   vanishes. Hence, it suffices to check that at any point $\mu\in
   \mathcal{O}_j\cap \mathcal{O}_{j'}\subset \mathfrak{t}\subset
   L\mathfrak{g}^*$, the contraction with $\zeta_{L\mathfrak{g}^*}$
   is equal for $\zeta\in L\mathfrak{g}$. Also, \begin{eqnarray*}
   \iota(\zeta)\sigma_j&=&\widetilde{\pi}_j^*dB(\Phi_j,\zeta)\\
   &=& dB(\widetilde{\pi}_j^*\Phi_j,\zeta)\end{eqnarray*}
   where $\Phi_j:\mathcal{O}_j\hookrightarrow L\mathfrak{g}^*$ is inclusion.
   But\begin{eqnarray*}(\widetilde{\pi}_j^*\Phi_j-\widetilde{\pi}_{j'}^*\Phi_{j'})
   \mid_{g\cdot\mu}&=&g\cdot\mu_j-g\cdot\mu_{j'}\\&=&(Ad_g(\mu_j)-g^*\theta^R)
   -(Ad_g(\mu_{j'})-g^*\theta^R)\\&=&
   Ad_g(\mu_j-\mu_{j'}).\end{eqnarray*}This, however, is another
   moment map for
   $\operatorname{Hol}^*(\varpi_j-\varpi_{j'})$.\end{proof}
   Thus, the locally defined forms
   $$\varpi\mid_{\widetilde{V}_j}=\operatorname{Hol}^*(\varpi_j)-\sigma_j$$
   patch together to define a global 2-form $\varpi\in
   \Omega^2(L\mathfrak{g}^*)$. Form the properties of $\varpi_j$ and
   $\sigma_j$ we read off:\\ \indent(i)
   $d\varpi=\operatorname{Hol}^*\eta$,\\ \indent(ii)$\iota(\zeta_{L\mathfrak{g}^*})\varpi=
   \frac{1}{2}B(\operatorname{Hol}^*(\theta^L+\theta^R),\zeta(0))-dB(
   \mu,\zeta)$. \\ Here, $\mu:L\mathfrak{g}^*\rightarrow
   L\mathfrak{g}^*$ is the identity map. Such a 2-form was
   constructed in~\cite{MR99k:58062} using a different method. \\ \indent
   Consider the case $G=SU(n)$ or $G=Sp(n)$, which are the two
   cases where the vertices of the alcove lie in the weight
   lattice $\Lambda^*$, where we identify $\mathfrak{g}^*\cong
   \mathfrak{g}$ using the basic inner product.
   Let$$U(1)\rightarrow \widehat{LG}\rightarrow LG$$ denote the
  $k$-th power of the basic central extension of the loop group~\cite{MR88i:22049}.
  That is, on the Lie algebra level the central
  extension$$\mathbb{R}\rightarrow\widehat{L\mathfrak{g}}\rightarrow
  L\mathfrak{g}$$ is defined by the cocycle,$$(\xi_1,\xi_2)\mapsto
  \int_{S^1}B(\xi_1,d\xi_2)\qquad \xi\in
  L\mathfrak{g}=\Omega^0(S^1,\mathfrak{g}))$$where $B$ is the
  inner product at level $k$. The coadjoint action of
  $\widehat{LG}$ on $\widehat{L\mathfrak{g}^*}=L\mathfrak{g}^*\times \mathbb{R}$
  preserves the level sets $L\mathfrak{g}^*\times\{t\}$,
  and the action for $t=1$ is exactly the gauge action of $LG$
  considered above.
   Since $\mu_j\in \Lambda^*$, the orbits $LG\cdot
   \mu_j=\mathcal{O}_j$ carry $\widehat{LG}$-equivariant
   pre-quantum line bundles
   $L_{\mathcal{O}_j}\rightarrow \mathcal{O}_j$, given explicitly
   as
   $$L_{\mathcal{O}_j}=\widehat{LG}\times_{(\widehat{LG})_j}\mathbb{C}_{(\mu_j,1)}.$$
   Here $(\widehat{LG})_j$ is the restriction of $\widehat{LG}$ to
   the stabilizer $(LG)_j$ of $\mu_j \in \mathfrak{t}\subset
   L\mathfrak{g}^*$, and $\mathbb{C}_{(\mu_j,1)}$ denotes the
   1-dimensional representation of $(\widehat{LG})_j$ with weight
   $(\mu_j,1)\in \Lambda^*\times \mathbb{Z}$. Let $E_j\rightarrow
   \widetilde{V}_j$ be the pull-back
   $\widetilde{\pi}_j^*L_{\mathcal{O}_j}$. On overlaps,
   $\widetilde{V}_j\cap \widetilde{V}_{j'}$, $E_j\otimes E_{j'}^{-1}$ is an
   associated bundle for the weight
   $(\mu_j,1)-(\mu_{j'},1)=(\mu_j-\mu_{j'},0)$. Therefore, $E_j\otimes
   E_{j'}^{-1}$ is an $LG$-equivariant bundle. It is clear by
   construction that $E_j\otimes E_{j'}^{-1}$ is the pull-back of the
   pre-quantum line bundle over $\mathcal{O}_{jj'}\subset
   \mathfrak{g}^*$. Taking all this information together, we have
   constructed an explicit quasi-line bundle for the pull-back of
   the $k$-th power of the basic gerbe under holonomy map
   $\operatorname{Hol}:L\mathfrak{g}^*\rightarrow G$ with error
   2-form equal to $\varpi\in \Omega^2(L\mathfrak{g}^*)$.
\section{Review of Quasi-Hamiltonian $G$-Spaces}
Suppose $(M,\omega)$ is a symplectic manifold together with a
symplectic action of a Lie group $G$. This action called Hamiltonian
if there exists a smooth equivariant map
\[\Phi: M\rightarrow \mathfrak{g}^*\]such that \[\iota
(\xi_M)\omega+d\langle\Phi,\xi\rangle=0\]for all $\xi\in
\mathfrak{g}$, where $\xi_M$ is the vector field on $M$ generated by
$\xi\in \mathfrak{g}$,
i.e.,\[\xi_M(m)=\frac{d}{dt}\mid_{t=0}\exp(t\xi)\cdot m.\] The map
$\Phi$ and the triple $(M,\omega,\Phi)$ are known as moment map and
Hamiltonian $G$-manifold respectively~\cite{MR91d:58073}. Let $G$ be
a compact Lie group. Fix an invariant inner product $B$ on
$\mathfrak{g}$, which we use to identify
$\mathfrak{g}\cong\mathfrak{g}^*$. Since the exponential map
$exp:\mathfrak{g}\rightarrow G$ is a diffeomorphism in a
neighborhood of the origin, the composition map \[\Psi:=exp\circ
\Phi:M\rightarrow G\] inherits the properties of the moment map
$\Phi$ and vice versa.
\begin{defn} A quasi-Hamiltonian $G$-space with group-valued moment map
is a triple $(M,\omega,\Psi)$ consisting of a $G$-manifold M, an
invariant 2-form $\omega\in\Omega^2(M)$, and an equivariant smooth
map $\Psi:M\rightarrow G$ such that
\begin{enumerate}\item $d\omega =\Psi^{*}\eta$ where $\eta \in
\Omega^{3}(G)$ is the 3-form defined by $B$. This condition is
called the relative cocycle condition.\item
$\iota(\xi_{M})\omega=\frac{1}{2}B(\Psi^{*}(\theta^{L}+\theta^{R}),\xi)$.
This condition is called the moment map condition.\item The $\ker(
\omega_{m})\in T_m(M)$ for $m\in M$ consists of all $\xi_{M}(m)$
such that
\[(Ad_{\Psi(m)}+1)\xi=0.\]This is called the minimal degeneracy
condition.\end{enumerate}
\end{defn}\subsection{Examples}
\begin{exmp}Consider a Hamiltonian $G$-manifold $(M,\omega,\Phi)$ such that
the image of $\Phi$ is a subset of the set of regular values for the
exponential map. Then $(M,\Upsilon,\Psi)$ is a quasi-Hamiltonian
$G$-space with group-valued moment map, where
\[\Psi=\exp\circ\Phi\]and
\[\Upsilon:=\omega+\Phi^*\varpi\]where $\varpi\in \Omega^2(\mathfrak{g})$
is the primitive for $\exp^*\eta$ given by the de Rham homotopy
operator for the vector space $\mathfrak{g}$. The converse is also
true, provided that $\Psi(M)$ lies in a neighborhood of the origin
on which the exponential map is a
diffeomorphism.\end{exmp}\begin{exmp}\label{Lisa} Let
$\mathcal{C}\subset G$ be a conjugacy class of $G$. The triple
$(\mathcal{C},\omega,\Phi)$ is a quasi-Hamiltonian $G$-space with
group-valued moment map where $\Phi:\mathcal{C}\hookrightarrow G$ is
inclusion and
$\omega_g(\xi_{\mathcal{C}}(g),\zeta_{\mathcal{C}}(g))=
\frac{1}{2}B((Ad_g-Ad_{g^{-1}})\xi,\zeta)$
~\cite{MR98e:58034}.\end{exmp}\begin{exmp}Given an involutive  Lie
group automorphism $\rho\in Aut(G)$, i.e., $\rho^2=1$, one defines
\emph{twisted conjugacy classes} to be the orbits of the action
$h\cdot g=\rho(h)gh^{-1}.$ $G$ is a symmetric space
$$G=G\times G/(G\times G)^\rho$$ where $\rho(g_1,g_2)=(g_2,g_1)$.
The map $G\times G\rightarrow \mathbb{Z}_2\ltimes G\times G,
(g_1,g_2)\mapsto(\rho^{-1},g_1,g_2)$ takes the twisted conjugacy
classes of $G\times G$ to conjugacy classes of the disconnected
group $\mathbb{Z}_2\ltimes G\times G$. Thus by using example
\ref{Lisa} the group $G$ itself becomes a group-valued Hamiltonian
$\mathbb{Z}_2\ltimes G\times G$, with 2-form $\omega=0$, moment map
$g\mapsto (\rho,g,g^{-1})$ and action $(g_1,g_2)\cdot g=g_2 g
g_1^{-1}, \rho\cdot g=g^{-1}$.\end{exmp}\begin{exmp}Let $D(G)$ be a
product of two copies of $G$. On $D(G)$, we can define a $G\times G$
action by\[(g_1,g_2).(a,b)=(g_1ag_2^{-1},g_2ag_1^{-1}).\]Define a
map\[\Psi:D(G)\rightarrow G\times G,\quad
\Psi(a,b)=(ab,a^{-1}b^{-1})\]and let the 2-form $\omega$ be defined
by
\[\omega=\frac{1}{2}(B(\operatorname{Pr}_1^{*}\theta^L,
\operatorname{Pr}_2^{*}\theta^R)+B(\operatorname{Pr}_1^{*}\theta^R,\operatorname{Pr}_2^{*}\theta^L))\]where
$\operatorname{Pr}_1$ and $\operatorname{Pr}_2$ are projections to
the first and second factor. Then the triple $(D(G),\omega,\Psi)$ is
a Hamiltonian $G\times G$-manifold with group-valued moment
map.\end{exmp}
\begin{exmp}\label{examplelisa}
Let $G=SU(2)$ and $M=S^4$ the unit sphere in
$\mathbb{R}^5\cong\mathbb{C}^2\times\mathbb{R}$, with SU(2)-action
induced from the action on $\mathbb{C}^2$. $M$ carries the structure
of a group-valued Hamiltonian $SU(2)$-manifold, with the moment map
$\Psi:M\rightarrow SU(2)\cong S^3$ the suspension of the Hopf
fiberation $S^3\rightarrow S^2.$ For details,
see~\cite{MR2003d:53151}. This example is generalized by
Hurtubise-Jeffrey-Sjamaar in~\cite{2} to $G=SU(n)$ acting on
$M=S^{2n}$ (viewed as unit sphere in $\mathbb{C}^n\times
\mathbb{R})$.
\end{exmp}
The equivariant de Rham complex is defined as
\[\Omega^{k}_G(M)= \bigoplus_{2i+j=k}(\,\Omega^{j}(M)\otimes S^i(\mathfrak{g}^*)
)^G\]where $S(\mathfrak{g}^*)$ is the symmetric algebra over the
dual of the Lie algebra of $G$. Elements in this complex can be
viewed as equivariant polynomial maps from $\mathfrak{g}$ into the
space of differential forms. $\Omega_G(M)$ carries an equivariant
differential $d_G$ of degree 1,
\[(d_G\alpha)(\xi):=d\alpha(\xi)+\iota(\xi_M)\alpha(\xi).\]Since
$(d+\iota(\xi_M))^2=L(\xi_M)$ and we are restricting on the
equivariant maps, $d_G^2=0$. The equivariant cohomology is the
cohomology of this co-chain complex~\cite{MR2001i:53140}. The
canonical 3-form $\eta$ has a closed equivariant extension
$\eta_G\in \Omega_G^3(G)$ given by
$$\eta_G(\xi):=\eta+\frac{1}{2}B(\theta^L+\theta^R,\xi).$$
We can combine the first two conditions of the definition of a
group-valued moment map and get the condition
$$d_G\omega=\Psi^*\eta_G.$$
\subsection{Products}Suppose $(M,\omega, (\Psi_1,\Psi_2))$ is a
group-valued Hamiltonian $G\times G$-manifold. Then
$\widetilde{M}=M$ with diagonal action, moment map
$\widetilde{\Psi}=\Psi_1\Psi_2$ and 2-form
$$\widetilde{\omega}=\omega-\frac{1}{2}
B(\Psi_1^{*}\theta^L,\Psi^*_{2}\theta^R)$$ is a group-valued
quasi-Hamiltonian $G$-space. If $\widetilde{M}=M_1\times M_2$ is a
direct product of two group-valued quasi-Hamiltonian $G$-spaces, we
call $\widetilde{M}$ the fusion product of $M_1$ and $M_2$. This
product is denoted by $M_1\circledast M_2$. If we apply fusion to
the double $D(G)$, we obtain a group-valued quasi-Hamiltonian
$G$-space with $G$-action
$$g\cdot(a,b)=(Ad_ga,Ad_gb),$$moment map
$$\Psi(a,b)=aba^{-1}b^{-1}\equiv[a,b],$$ and 2-form
$$\omega=\frac{1}{2}(B(\operatorname{Pr}_1^{*}\theta^L,\operatorname{Pr}_2^{*}\theta^R)
+B(\operatorname{Pr}_1^{*}\theta^R,\operatorname{Pr}_2^{*}\theta^L)-
B((ab)^{*}\theta^L,(a^{-1}b^{-1})^*\theta^R)).$$Fusion of $h$ copies
of $D(G)$ and conjugacy classes $\mathcal{C}_1,\cdots,\mathcal{C}_r$
gives a new quasi-Hamiltonian space with the moment
map\[\Psi(a_1,b_1,\cdots,a_h,b_h,d_1,\cdots,d_r)=\prod^h_{j=1}[a_j,b_j]\prod^r_{k=1}d_k.\]
\subsection{Reduction}The symplectic reduction works as usual:
\\ \indent If $(M,\omega,\Psi)$ be a Hamiltonian $G$-space with group-valued moment map
and the identity element $e\in G$ be a regular value of $\Psi$, then
$G$ acts locally freely on $\Psi^{-1}(e)$ and therefore
$\Psi^{-1}(e)/G$ is smooth. Furthermore, the pull- back of $\omega$
to identity level set descends to a symplectic form on
$M//G:=\Psi^{-1}(e)/G$. For instance, the moduli space of flat
G-bundles on a closed oriented surface of genus h with $r$ boundary
components, can be written
\[\mathcal{M}(\Sigma;\mathcal{C}_1,\cdots,\mathcal{C}_r)=G^{2h}\circledast\mathcal{C}_1
\circledast\cdots\circledast\mathcal{C}_r//G= \Psi^{-1}(e)/G\]where
the $j$-th boundary component is the bundle corresponding to the
conjugacy class $\mathcal{C}_j$. More details can be found in
~\cite{MR99k:58062},~\cite{MR2003d:53151}.

\section{Pre-quantization of Quasi-Hamiltonian $G$-Spaces}
We know that a symplectic manifold $(M,\omega)$ is pre-quantizable
(admits a line bundle L over M with curvature 2-form $\omega$) if
the 2-form $\omega$ is integral. In this Section, we will first
introduce a notion of a pre-quantization of a space with G-valued
moment map and then  give a similar criterion for being
pre-quantizable.\begin{defn} Let G be a compact connected Lie group
with canonical 3-form $\eta$. Fix a gerbe $\mathcal{G}$ on $G$ with
connection $(\nabla,\varpi)$ such that $curv(\mathcal{G})=\eta$. A
pre-quantization of $(M,\omega,\Psi)$ is a relative gerbe with
connection $(\mathcal{L},\g)$ corresponding to the map $\Psi$ with
relative curvature $(\omega,\eta)$.
\end{defn}Since $\eta$ is closed 3-form and $\Psi^*\eta=d\omega$,
$(\omega,\eta)$ defines a relative cocycle. Recall from chapter 1
that a class $[(\omega,\eta)]\in H^3(\Psi,\mathbb{R})$ is integral
if and only if $\int_{\beta}\eta-\int_\Sigma\omega\in \mathbb{Z}$
for all relative cycles $(\beta,\Sigma)\in C_3(\Psi,\mathbb{R})$.
\begin{rem}
\label{integral} $(M,\omega,\Psi)$ is pre-quantizable if and only if
$[(\omega,\eta)]$ is integral by Theorem \ref{quantization}.
\end{rem}
\begin{thm}Suppose $M_i$, $i=1,2$ are two
quasi-Hamiltonian $G$-spaces. The fusion product $M_1\circledast
M_2$ is pre-quantizable if both $M_1$ and $M_2$ are
pre-quantizable.\end{thm}\begin{proof}Let $\operatorname{Mult}:
G\times G\rightarrow G$ be group multiplication and
$\operatorname{Pr}_i: G\times G\rightarrow G$, $i=1,2$ projections
to the first and second factors.
   Since $$\operatorname{Mult}^*\eta=\operatorname{Pr}_1^*\eta+\operatorname{Pr}_2^*\eta
   +\frac{1}{2}B(\operatorname{Pr}_1^*\theta^L,\operatorname{Pr}_1^*\theta^R),$$
   we get a quasi-line bundle with connection for the gerbe
   $\operatorname{Mult}^*\g\otimes (\operatorname{Pr}_1^*\g)^{-1}\otimes
   (\operatorname{Pr}_2^*\g)^{-1}$ such that the error 2-form is
   equal to
   $\frac{1}{2}B(\operatorname{Pr}_1^*\theta^L,\operatorname{Pr}_1^*\theta^R)$.
   Any two such quasi-line bundles differ by a flat line bundle
   with connection. Let $\Psi_i$, $i=1,2$ be moment maps
   for $M_i$, $i=1,2$ respectively and $\Psi=\Psi_1\Psi_2$
   be the moment map for their fusion product $M_1\circledast
M_2$. Thus,\begin{eqnarray*}\Psi^*\g &=&
   (\Psi_1\times \Psi_2)^*\operatorname{Mult}^*\g\\&=& (\Psi_1^*\times \Psi_2^*)
   \big((\operatorname{Pr}_1^*\g)\otimes
   (\operatorname{Pr}_2^*\g)\big)\\&=&\Psi_1^*\g\otimes
   \Psi_2^*\g.\end{eqnarray*}Therefore $M_1\circledast
M_2$ is pre-quantizable if and only if both $M_1$ and $M_2$ are
pre-quantizable.\end{proof}
\begin{prop}Suppose $G$ is simple and simply connected. Let
$k\in\mathbb{Z}$ be the level of $(M,\omega,\Psi)$. Suppose
$H^2(M,\mathbb{Z})=0$. Then there exists a pre-quantization of
$(M,\omega,\Psi)$ if and only if the image of
$$\Psi^*:H^3(G,\mathbb{Z})\rightarrow
H^3(M,\mathbb{Z})$$ is $k$-torsion.\end{prop}
\begin{proof}By assumption, $[\eta]$ represents $k$ times the generator of
$H^3(G,\mathbb{Z})$. If $H^2(M,\mathbb{Z})=0$, the long exact
sequence:\[\cdots\rightarrow H^2(G,\mathbb{Z})\rightarrow
H^2(M,\mathbb{Z})\rightarrow H^3(\Psi,\mathbb{Z})\rightarrow
H^3(G,\mathbb{Z})\overset{\Psi^*}\rightarrow
H^3(M,\mathbb{Z})\rightarrow \cdots\]shows that the map
$H^3(\Psi,\mathbb{Z})\rightarrow H^3(G,\mathbb{Z})$ is injective. In
particular, $H^3(\Psi,\mathbb{Z})$ has no torsion, and
$(M,\omega,\Psi)$ is pre-quantizable if and only if $[\eta]$ is in
the image of the map $H^3(\Psi,\mathbb{Z})\rightarrow
H^3(G,\mathbb{Z})$, i.e., in the kernel of
$H^3(G,\mathbb{Z})\rightarrow H^3(M,\mathbb{Z})$. This exactly means
that the image of this map is $k$-torsion.\end{proof}
\begin{prop}\label{condition}
If $\,H_2(M,\mathbb{Z})=0$ a pre-quantization exists. More
generally, if $H_2(M,\mathbb{Z})$ is r-torsion, a level $k$
pre-quantization exists, where $k$ is a multiple of $r$.\end{prop}
\begin{proof}If $rH_2(M,\mathbb{Z})=0$, for any
cycle $S\in C_2(M)$, there is a 3-chain $T\in C_3(M)$ with $\partial
T=r\cdot S$. If $\Psi(S)=\partial B$, $\Psi(T)-rB$ is a cycle and
\begin{equation}\begin{split}\int_S
k\omega -\int_B k\eta &=\frac{k}{r}(\int_Td\omega-\int_{rB}\eta)\\
&=\frac{k}{r}(\int_T \Psi^*\eta-\int_{rB}\eta)\\
&=\frac{k}{r}(\int_{\Psi (T)}\eta-\int_{rB}\eta)\\
&=\frac{k}{r}(\int_{\Psi(T)-rB}\eta)\in
\mathbb{Z}.\end{split}\end{equation}By Remark \ref{integral}
$(M,\omega,\Psi)$ is
pre-quantizable.\end{proof}\begin{exmp}$M=S^{4}$ carries the
structure of a group-valued Hamiltonian $SU(2)$-manifold, with the
moment map $\Psi:M\rightarrow SU(2)\cong S^3$ the suspension of the
Hopf fiberation $S^3\rightarrow S^2$, Example \ref{examplelisa}. By
Proposition \ref{condition} this $SU(2)$-valued moment map is
pre-quantizable.\end{exmp}
\subsection{Reduction}Let $G$ be a simply connected Lie group. Fix
a pre-quantization $\lc$ for a space with G-valued moment map
$(M,\omega,\Psi)$.\[\\[7.pt]\begin{CD}
\lc  @. \g\\
@VVV   @VVV\\
M @>\Psi>>G\\
@A{\iota}AA   @AAA\\
\Psi^{-1}(e) @>\Psi>> \{e\}
\end{CD}\\[7.pt]\]
Since $\mathcal{G}|_{\Psi^{-1}(e)}$ is equal to trivial gerbe,
$\mathcal{L}|_{\Psi^{-1}(e)}$ is a line bundle with connection with
curvature $(\iota_{\Psi^{-1}(e)})^*\omega$. Since G is simply
connected and the 2-form $(\iota_{\Psi^{-1}(e)})^*\omega$ is
$G$-basic, there exists a unique lift of the $G$-action to
$\mathcal{L}|_{\Psi^{-1}(e)}$ in such a way that the generating
vector fields on $\mathcal{L}|_{\Psi^{-1}(e)}$ are horizantal. This
is a special case of Kostant's construction~\cite{MR45:3638}. In
conclusion, we get a pre-quantum line bundle over $\Psi^{-1}(e)/G$.
\subsection{A Finite Dimensional Pre-quantum Line Bundle for
$\mathcal{M}(\Sigma)$} Let $M=G^{2h}$ where $G$ is a simply
connected Lie group and consider the map
$$\Psi:M\rightarrow G$$with the rule
\[\Psi(a_1,\cdots,a_h)= \prod_{i=1}^h[a_i,b_i].\]Let $\mathcal{G}$
be the basic gerbe with the connection on G and $curv({\g})=\eta$.
The moduli space of flat G-bundles on a closed oriented surface
$\Sigma$ of genus h is equal to
\[\mathcal{M}(\Sigma)=G^{2h}//G=\Psi^{-1}(e)/G.\]Since $G$ is simply
connected, $H^2(G,\mathbb{Z})=H^2(G^{2h},\mathbb{Z})=0$.
$H^3(G^{2h},\mathbb{Z})\simeq \mathbb{Z}$ is torsion free therefore
by Proposition 4.3.3 there exists a unique quasi-line
bundle$\mathcal{L}$ for the gerbe $\Psi^*\mathcal{G}$.
\newline \indent Pick a connection for this quasi-line bundleand call the error 2-form
$\nu$. Therefore $d(\nu-\omega)=0$. This together with the fact that
$H^2(M,\mathbb{Z})=0$ allow us to modify quasi-line bundle with
connection such that triple $(M,\omega,\Psi)$ is pre-quantizabe. By
reduction we get a pre-quantum line bundle over $\Psi^{-1}(e)/G=
\mathcal{M}(\Sigma)$.

\subsection{Pre-quantization of Conjugacy Classes of a Lie
Group}Let G be a simple, simply connected compact Lie group. Fix an
inner product $B$ at level $k$. The map
$$\exp:\mathfrak{g}\rightarrow G$$ takes
(co)adjoint orbits $\mathcal{O}_{\xi}$ to conjugacy classes
$\mathcal{C}=G\cdot{\exp(\xi)}$.\\ \indent Any conjugacy class
$\mathcal{C}\subseteq G$ is uniquely a G-valued quasi-Hamiltonian
$G$-space $(\,\mathcal{C},\omega,\Psi)$, where
$\Psi:\mathcal{C}\hookrightarrow G$ is inclusion map, as it
explained in example \ref{Lisa}. Suppose $(\beta,\Sigma)\in
\Co_n(\Psi,\mathbb{Z})$ is a cycle. We want to see under which
conditions $\mathcal{C}$ is pre-quantizable at level $k$.
Equivalently, we are looking for conditions which implies
\[k(\int_{\beta}\eta-\int_\Sigma\omega)\in \mathbb{Z}\] where
$\eta=\frac{1}{12}B(\theta^L,[\theta^L,\theta^L])$ is canonical
3-form. Consider the basic gerbe
$\mathcal{G}=(\mathcal{V},L,\theta)$ with connection on G with
curvature $\eta$. For all $\mathcal{C}\subseteq G$ there exists a
unique $\xi\in\mathfrak{A}$ such that $\exp(\xi)\in \mathcal{C}$.
Let
$$\iota_{\mathcal{C}}:\mathcal{C}\rightarrow G$$ be inclusion map
assume that $\varpi_0$ is the primitive of $\eta$ on $V_0$, i.e.,
$$\eta\mid_{V_0}=d\varpi_{0}$$where $V_{0}$ contains $\mathcal{C}$.
Recall from Section \ref{basic gerbe} that
$$\omega_{\mathcal{C}}=\theta_0^*(\omega_{\mathcal{O}_{\xi}})_G-\iota_{\mathcal{C}}^*(\varpi_0)_G$$
and pull-back of the $\theta$ to $\mathcal{C}$ is zero. Thus,
\begin{eqnarray*}k(\int_{\Sigma}\omega_{\mathcal{C}}-\int_{\beta}\iota_{\mathcal{C}}^*\eta) &=&
k(\int_{\Sigma}\theta_0^*(\omega_{\mathcal{O}_{\xi}})_G-\iota_{\mathcal{C}}^*(\varpi_0)_G-\int_{\beta}
\iota_{\mathcal{C}}^*(\eta))\\
&=&k(\int_{\Sigma}\theta_0^*(\omega_{\mathcal{O}_{\xi}})_G).\end{eqnarray*}
$(\omega_{\mathcal{C}},\iota_{\mathcal{C}}^*\eta)$ is integral if
and only if the symplectic 2-form
 $k\omega_{\mathcal{O}_{\xi}}$ is
integral. It is a well-known fact from symplectic geometry that
$k\omega_{\mathcal{O}_{\xi}}$ is integral if and only if $B(\xi)\in
\Lambda^*_k:=\Lambda^*\cap k\mathfrak{A}$, by viewing $B$ as a
linear map $\mathfrak{t}\rightarrow \mathfrak{t}^*$.
   \section{Hamiltonian Loop Group Spaces}
   Fix an invariant inner product $B$ on $\mathfrak{g}$.
   Assume that $G$ is simple and simply connected.
   Recall that a Hamiltonian loop group manifold is a triple
   $(\widetilde{M},\widetilde{\omega},\widetilde{\Psi})$ where
   $\widetilde{M}$ is an (infinite-dimensional) $LG$-manifold,
   $\widetilde{\omega}$ is an invariant symplectic form on $\widetilde{M}$,
   and $\widetilde{\Psi}:\widetilde{M}\rightarrow L\mathfrak{g}^*$
   an equivariant map satisfying the usual moment map condition,
   $$\iota(\xi_{\widetilde{M}})\widetilde{\omega}+dB(\widetilde{\Psi},\xi)=0
   \qquad\qquad\xi\in\Omega^0(S^1,\mathfrak{g}).$$
   \begin{exmp}Let $\mathcal{O}\subset
   L\mathfrak{g}^*$
   be an orbit of the loop group action. Then $\mathcal{O}$
   carries a unique structure for a Hamiltonian $LG$-manifold when
   the moment map is inclusion and the 2-form is $$
   \widetilde{\omega}_{\mu}(\xi_{\mathcal{O}}(\mu),\eta_{\mathcal{O}}(\mu))=\langle
   d_{\mu}\xi,\eta\rangle=\int_{S^1}B((d_{\mu}\xi),\eta).$$\end{exmp}
   The based loop group $\Omega G\subset LG$ consisting of loops that are trivial
   at the origin of $S^1$, acts freely on $L\mathfrak{g}^*$
   and the quotient is just the holonomy map. There
   is a one-to-one correspondence between quasi-Hamiltonian $G$-spaces
   $(M,\omega,\Psi)$ and Hamiltonian $LG$-spaces with proper
   moment maps
   $(\widetilde{M},\widetilde{\omega},\widetilde{\Psi})$, where
   \begin{eqnarray*}M &=& \widetilde{M}/\Omega G,\\
   \operatorname{Hol}\circ \widetilde{\Psi}&=& \Psi\circ
   \operatorname{Hol},\\
   \widetilde{\omega}&=&\operatorname{Hol}^*\omega-\widetilde{\Psi}\varpi.
   \end{eqnarray*}This is called \emph{Equivalence Theorem} in~\cite{MR99k:58062}.
   We thus, have a commutative diagram:
   \\*\[\begin{CD}
   \widetilde{M}  @>\widetilde{\Psi}>>     L\mathfrak{g}^*\\
   @V{\operatorname{Hol}}VV   @V{\operatorname{Hol}}VV\\
   M@>\Psi>>G
   \end{CD}\]\\\begin{thm}(Equivalence Theorem for Pre-quantization) There is a one-to-one correspondence between
   pre-quantizations of quasi-Hamiltonian $G$-spaces with group
   valued moment maps and pre-quantizations of the corresponding Hamiltonian $LG$-spaces
with proper moment maps.\end{thm}\begin{proof} Assume that we have
constructed a pre-quantization of a quasi-Hamiltonian $G$-space with
group-valued moment map $(M,\omega,\Psi)$ with the corresponding
Hamiltonian $LG$-space
$(\widetilde{M},\widetilde{\omega},\widetilde{\Psi})$. Thus, we have
a relative gerbe mapping to the basic gerbe over $G$. Pull-back of
this quasi-line bundle under
$\operatorname{Hol}:\widetilde{M}\rightarrow M$, gives a quasi-line
bundle of the gerbe
$\operatorname{Hol}^*\Psi^*\g=\widetilde{\Psi}^*\operatorname{Hol}^*\g$
over $\widetilde{M}$. But recall that, we have a quasi-line bundle
for $\operatorname{Hol}^*\g$ as it explained in Section \ref{rhol}.
Therefore the difference between these two quasi-line bundles with
connection is a line bundle with connection
$\widetilde{L}\rightarrow \widetilde{M}$ with the curvature 2-form
$\widetilde{\omega}=\operatorname{Hol}^*\omega-\widetilde{\Psi}\varpi$
by Remark \ref{qc}. Note also that if the quasi-line bundle for
$\Psi:M\rightarrow G$ is $G$-equivariant, then since the quasi-line
bundle for $L\mathfrak{g}^*\rightarrow G$ is
$\widehat{LG}$-equivariant, the line bundle $\widetilde{L}$ will be
$\widehat{LG}$ equivariant. Conversely, suppose that we are given a
$\widehat{LG}$-equivariant line bundle over $\widetilde{M}$, where
$U(1)\subset \widehat{LG}$ acts with weight 1. The difference of
this $\widehat{LG}$-equivariant line bundle and the
$\widehat{LG}$-equivariant quasi-line bundle for
$\operatorname{Hol}^*\Psi^*\g$ (constructed in Section \ref{rhol}),
is a quasi-line bundle with error 2-form
$\operatorname{Hol}^*\omega$. By descending of this quasi-line
bundle to $M$, we can get the desired quasi-line bundle for
$\Psi^*\g$.
\end{proof}
The argument, given here applies in greater generality:\newline For
any $\widehat{LG}$-equivariant line bundle $\widetilde{L}\rightarrow
\widetilde{M}$, where the central extension
$U(1)\subset\widehat{LG}$ acts with weight $k\in\mathbb{Z}$, there
is a corresponding relative gerbe at level $k$ with respect to the
map $\Psi:M\rightarrow G$. Indeed, the given quasi-line bundle for
$\widetilde{\Psi}^*\operatorname{Hol}^*\g^k$ is given by
$\widehat{LG}$ equivariant line bundles over
$\operatorname{Hol}^{-1}\Psi^{-1}(V_j)$ at level $k$. Twisting by
$\widetilde{L}$, we get new quasi-line bundle where $U(1)\subset
\widehat{LG}$ acts trivially. The quotient therefore descends to a
quasi-line bundle over $M$. For instance, Meinrenken and Woodward
construct for any Hamiltonian loop group space a so-called
``canonical line bundle'' in~\cite{MR2002g:53151}, which is
$\widehat{LG}$-equivariant at level $2c$, where $c$ is the dual
Coxeter number. Therefore this line bundle gives rise to a
distinguished element of $H^3(\Phi,\mathbb{Z})$ at level $2c$.
Notice that $M$ and $\widetilde{M}$ are not pre-quantizable
necessarily.


\subsection*{Acknowledgment}
The results of this paper were obtained during the author's Ph.D.
study at the University of Toronto and were also included in her
thesis dissertation. The guidance and support received from Eckhard
Meinrenken are deeply acknowledged.

\bibliographystyle{amsplain}
\bibliography{xbib}

\providecommand{\bysame}{\leavevmode\hbox to3em{\hrulefill}\thinspace}
\providecommand{\MR}{\relax\ifhmode\unskip\space\fi MR }
\providecommand{\MRhref}[2]{%
  \href{http://www.ams.org/mathscinet-getitem?mr=#1}{#2}
}
\providecommand{\href}[2]{#2}
\begin{thebibliography}{10}

\bibitem{MR2003d:53151}
A.~Alekseev, E.~Meinrenken, and C.~Woodward, \emph{Duistermaat-{H}eckman
  measures and moduli spaces of flat bundles over surfaces}, Geom. Funct. Anal.
  \textbf{12} (2002), no.~1, 1--31. \MR{2003d:53151}

\bibitem{MR99k:58062}
Anton Alekseev, Anton Malkin, and Eckhard Meinrenken, \emph{Lie group valued
  moment maps}, J. Differential Geom. \textbf{48} (1998), no.~3, 445--495.
  \MR{99k:58062}

\bibitem{MR1968268}
Kai Behrend, Ping Xu, and Bin Zhang, \emph{Equivariant gerbes over compact
  simple {L}ie groups}, C. R. Math. Acad. Sci. Paris \textbf{336} (2003),
  no.~3, 251--256. \MR{1 968 268}

\bibitem{MR39:1590}
N.~Bourbaki, \emph{\'{E}l\'ements de math\'ematique. {F}asc. {XXXIV}. {G}roupes
  et alg\`ebres de {L}ie. {C}hapitre {IV}: {G}roupes de {C}oxeter et syst\`emes
  de {T}its. {C}hapitre {V}: {G}roupes engendr\'es par des r\'eflexions.
  {C}hapitre {VI}: syst\`emes de racines}, Actualit\'es Scientifiques et
  Industrielles, No. 1337, Hermann, Paris, 1968. \MR{39 \#1590}

\bibitem{D}
D.~Chatterjee, \emph{On the construction of abelian gerbs}, Ph.D. thesis,
  University of Cambridge, 1998.

\bibitem{DK}
J.~J. Duistemaat and J.~A.~C. Kolk, \emph{Lie groups}, Springer-Verlag, Berlin,
  2000.

\bibitem{MR2003m:81222}
Krzysztof Gawedzki and Nuno Reis, \emph{W{ZW} branes and gerbes}, Rev. Math.
  Phys. \textbf{14} (2002), no.~12, 1281--1334. \MR{2003m:81222}

\bibitem{MR91d:58073}
Victor Guillemin and Shlomo Sternberg, \emph{Symplectic techniques in physics},
  second ed., Cambridge University Press, Cambridge, 1990. \MR{91d:58073}

\bibitem{MR2001i:53140}
Victor~W. Guillemin and Shlomo Sternberg, \emph{Suppersymmetry and equiavriant
  de {R}ham theory}, Mathematics Past and Present, Springer-Verlag, Berlin,
  1999, With an appendix containg two reprints by Henri Cartan [MR {\bf
  13},107e; MR {\bf 13},107f. \MR{2001i:53140}

\bibitem{MR98e:58034}
K.~Guruprasad, J.~Huebschmann, L.~Jeffrey, and A.~Weinstein, \emph{Group
  systems, groupoids and moduli spaces of parabolic bundles}, Duke Math. J.
  \textbf{89} (1997), no.~2, 377--412. \MR{98e:58034}

\bibitem{1}
N.~Hitchin, \emph{What is a gerbe?}, Notices of the A.M.S (2003), 218--219.

\bibitem{MR2003f:53086}
Nigel Hitchin, \emph{Lectures on special {L}agrangian submanifolds}, Winter
  School on Mirror Symmetry, Vector Bundles and Lagrangian Submanifolds
  (Cambridge, MA, 1999), AMS/IP Stud. Adv. Math., vol.~23, Amer. Math. Soc.,
  Providence, RI, 2001, pp.~151--182. \MR{2003f:53086}

\bibitem{2}
J~Hurtubise, L.~Jeffrey, and R.~Sjamaar, \emph{Group-valued implosion snd
  parabolic structure \textrm{I}}, math.S6/0402464.

\bibitem{MR2003c:22001}
Anthony~W. Knapp, \emph{Lie groups beyond an introduction}, second ed.,
  Progress in Mathematics, vol. 140.

\bibitem{MR45:3638}
Bertram Kostant, \emph{Quantization and unitary representations. {I}.
  {P}requantization},  (1970), 87--208. Lecture Notes in Math., Vol. 170.
  \MR{45 \#3638}

\bibitem{MR2002g:53151}
E.~Meinrenken and C.~Woodward, \emph{Canonical bundles for {H}amiltonian loop
  group manifolds}, Pacific J. Math. \textbf{198} (2001), no.~2, 477--487.
  \MR{2002g:53151}

\bibitem{EM}
Eckhard Meinrenken, \emph{The basic gerbe over a compact simple lie group},
  L'Enseignement Mathematique \textbf{49} (2003), 307--333.

\bibitem{MR88i:22049}
Andrew Pressley and Graeme Segal, \emph{Loop groups}, Oxford Mathematical
  Monographs, The Clarendon Press Oxford University Press, New York, 1986,
  Oxford Science Publications. \MR{88i:22049}

\bibitem{ZI}
Z.~Shahbazi, \emph{Relative gerbes}, To appear in Geometry and Physics
  (math.SG/0508016) (2005).

\end{thebibliography}
\end{document}